\documentclass[preprint]{elsarticle}

\usepackage{graphicx}
\usepackage{epstopdf}

\usepackage{amssymb}
\usepackage{amsthm}
\usepackage{amsmath}
\usepackage{url}
\usepackage{hyperref}
\usepackage{multirow}
\usepackage{epstopdf}
\usepackage{booktabs}
\usepackage{mathrsfs}

\journal{Elsevier}

\begin{document}

\begin{frontmatter}



\title{Numerical solution of space-fractional partial differential equations by a differential quadrature approach}

\cortext[cor1]{Corresponding author}
\author{X. G. Zhu}
\author{Y. F. Nie\corref{cor1}}
\ead{yfnie@nwpu.edu.cn}

\address{Department of Applied Mathematics, Northwestern Polytechnical University, Xi' an 710129, P.R. China}

\begin{abstract}
This article aims to develop a direct numerical approach to solve the space-fractional partial differential equations (PDEs) based on
a new differential quadrature (DQ) technique. The fractional derivatives are approximated by 
the weighted linear combinations of the function values at discrete grid points on problem domain with the weights  
calculated via using three types of radial basis functions (RBFs) as test functions.   
The method in presence is robust, straight forward to apply, and highly accurate under the condition that
the shape parameters of RBFs are well chosen. Numerical tests are provided to illustrate its validity and capability.
\end{abstract}


\end{frontmatter}


\newtheorem{theorem}{Theorem}[section]
\newtheorem{lemma}{Lemma}[section]
\newtheorem{definition}{Defition}[section]
\newtheorem{remark}{Remarks}[section]
\newtheorem{assu}{Assumption}[section]
\renewcommand{\theequation}{\arabic{section}.\arabic{equation}}

\section{Introduction}\label{s1}
In this study, we are mainly interested in an efficient method for numerically
solving a class of space-fractional models in the following form  
\begin{align}
\frac{\partial y(x,t)}{\partial t}-\kappa(x)\frac{\partial^\alpha y(x,t)}{\partial_+ x^\alpha}-
\upsilon(x)\frac{\partial^\alpha y(x,t)}{\partial_- x^\alpha} =f(x,t), \ \ x\in\varLambda, \ \ 0<t\leqslant T,\label{eq01}
\end{align}
subjected to the initial and boundary conditions  
\begin{align}
  & y(x,0)=\psi(x),\quad x\in\varLambda, \label{eq02}\\
  & y(a,t)=g_1(t), \quad y(b,t)=g_2(t), \quad 0<t\leqslant T, \label{eq03}
\end{align} 
where $1<\alpha\leqslant2$, $\varLambda=[a,b]$, $\kappa(x)$, $\upsilon(x)$ are non-negative but do not vanish altogether.
$g_1(t)\neq0$ only when $\kappa(x)\equiv0$ and $g_2(t)\neq0$ only when $\upsilon(x)\equiv0$. In Eq. (\ref{eq01}), the
space-fractional derivatives are defined in Caputo sense, i.e.,
\begin{align*}
 \frac{\partial^\alpha y(x,t)}{\partial_+ x^\alpha}&=\frac{1}{\Gamma(2-\alpha)}
    \int^x_a\frac{\partial^2 y(\xi,t)}{\partial \xi^2}\frac{d\xi}{(x-\xi)^{\alpha-1}},\\
 \frac{\partial^\alpha y(x,t)}{\partial_- x^\alpha}&=\frac{1}{\Gamma(2-\alpha)}
    \int^b_x\frac{\partial^2 y(\xi,t)}{\partial \xi^2}\frac{d\xi}{(\xi-x)^{\alpha-1}},
\end{align*}
with the Euler's Gamma function $\Gamma(\cdot)$.  

The space-fractional PDEs describe many physical phenomena  
such as anomalous transport, hereditary elasticity, and chaotic dynamics \cite{R01,R02,R03}, while compared favorably to the integer PDEs.
Since they are frequently sufficing in the absence of exact closed-form solutions,
various numerical algorithms have been designed to solve them, typically including 
general Pad\'{e} approximation \cite{R04}, finite difference methods \cite{R05,R07,R06}, 
meshless point interpolation method \cite{R09}, finite element methods \cite{R12,R13}, 
discontinuous Galerkin method \cite{R14}, finite volume method \cite{R15}, 
spline approximation method (SAM) \cite{R17}, 
RBF Kansa method \cite{R11}, polynomial and fractional spectral collocation methods \cite{R19,R20}.
In \cite{R24,R21,R23,R22}, a series of operational matrix methods are constructed
via the approximate expansions using shifted Jacobi, Chebyshev, Legendre polynomials, 
and Haar wavelets functions, as elements, respectively. Some analytic techniques are referred to \cite{R27,R25,R28,R26}
and references therein.

DQ method is understood as a direct numerical approach for PDEs that evaluates the derivatives via representative weighted
linear combinations of function values on problem domain \cite{R30}. A group of test functions to calculate these weights
can be chosen as Lagrange basis functions,  RBFs, and orthogonal polynomials \cite{R60,R37,R39,R36}. DQ method is characterized
by a few advantages such as high accuracy,\ low occupancy cost,\ truly \emph{mesh-free} and the ease of programming.

Due to the non-locality of fractional derivatives, a great extra computational cost is usually incurred
when a conventional algorithm is applied to a fractional PDE.
In this work, we propose a new RBFs based DQ method (RBF-DQM) for Eqs. (\ref{eq01})-(\ref{eq03}). 
Using three types of RBFs as test functions, the weights are successfully determined and with them, the
equation under consideration degenerates to an ordinary differential system (ODS).
A time-stepping RBF-DQM is derived by introducing a difference scheme in time.
The presented technique inherits the features of classic DQ methods. More importantly, it is insensitive to dimensional change, so
it serves as a good alternative for the high-dimensional or the other complex fractional models arising in actual applications.

The outline is as follows. In Section \ref{s2}, we give a brief description of fractional derivatives. 
In Section \ref{s3}, the weighted coefficients are calculated by commonly used RBFs,
which are required to approximate the fractional derivatives. We propose a Crank-Nicolson RBF-DQM to discretize
the model problem in Section \ref{s4} and test its codes on three illustrative examples in Section \ref{s5}.
A conclusion is drawn in the last section.

\section{Fractional derivatives}\label{s2}
At first, some basic definitions are introduced for preliminaries.
Let $\alpha\in \mathbb{R}^+$; then the following formulas
\begin{align*}
   {_a}D_x^\alpha f(x)&=\frac{1}{\Gamma(m-\alpha)}
    \int^x_a\frac{\partial^m f(\xi)}{\partial \xi^m}\frac{d\xi}{(x-\xi)^{\alpha-m+1}},\\
   {_x^*}D_b^\alpha f(x)&=\frac{(-1)^m}{\Gamma(m-\alpha)}
    \int^b_x\frac{\partial^m f(\xi)}{\partial \xi^m}\frac{d\xi}{(\xi-x)^{\alpha-m+1}},
\end{align*}
define the left and right $\alpha$-th Caputo derivatives, respectively, if
$f(x)\in C^{m}(\varLambda)$, where $m=[\alpha]+1$ for $\alpha\notin\mathbb{N}$, $m=\alpha$
for $\alpha\in\mathbb{N}$, and $[\cdot]$ is the floor function.

The left and right Caputo derivatives have the properties
\begin{align*}
   {_a}D_x^\alpha (x-a)^\beta=\frac{\Gamma(\beta+1)}{\Gamma(\beta-\alpha+1)}(x-a)^{\beta},\quad
   {_x^*}D_b^\alpha (b-x)^\beta=\frac{\Gamma(\beta+1)}{\Gamma(\beta-\alpha+1)}(b-x)^{\beta},
\end{align*}
and coincide with classic derivatives with exactness to a multiplier factor $(-1)^s$:
\begin{align*}
   {_0}D_x^s f(x)=\frac{\partial^s f(x)}{\partial x^s},\quad
   {_x^*}D_b^s f(x)=(-1)^s\frac{\partial^s f(x)}{\partial x^s}, 
\end{align*}
where $\beta>m-1$ and $s\in \mathbb{N}$. We refer the readers to \cite{R29,R35} for more properties.

\section{DQ formulations based on RBFs}\label{s3}
In the sequel, DQ formulations for fractional derivatives based on RBFs are derived.
Define a lattice on $[a,b]$ but not necessarily with equally spaced points, i.e.,
$a=x_0<x_1<\cdots< x_{M-1}<x_M=b$, $M\in \mathbb{Z}^+$. In general, we always approximate the exact solution
of a PDE like Eqs. (\ref{eq01})-(\ref{eq03}) in the form 
\begin{equation}\label{eq04}
    y(x,t)\cong\sum_{k=0}^M\delta_k(t)\phi_k(x), 
\end{equation}
with a set of proper basis functions $\{\phi_k(x)\}_{k=0}^M$. However, if
\begin{align}
\frac{\partial^\alpha \phi_k(x_i)}{\partial_+ x^\alpha}&=\sum\limits_{j=0}^M {a_{ij}^{(\alpha)}\phi_k(x_j)},\ \ i,k=0,1,\ldots M,\label{eq05}\\
\frac{\partial^\alpha \phi_k(x_i)}{\partial_- x^\alpha}&=\sum\limits_{j=0}^M {b_{ij}^{(\alpha)}\phi_k(x_j)},\ \ i,k=0,1,\ldots M,\label{eq06}
\end{align}
and on acting 
$\frac{\partial^\alpha}{\partial_+ x^\alpha}$, $\frac{\partial^\alpha}{\partial_- x^\alpha}$ on both sides of Eq. (\ref{eq04}), we realize that
\begin{small}
\begin{align}
\frac{\partial^\alpha y(x_i,t)}{\partial_+ x^\alpha}\cong\sum_{k=0}^M\delta_k(t)\frac{\partial^\alpha \phi_k(x_i)}{\partial_+ x^\alpha}
        =\sum_{k=0}^M\delta_k(t)\sum\limits_{j=0}^M {a_{ij}^{(\alpha)}\phi_k(x_j)}\cong \sum\limits_{j=0}^M {a_{ij}^{(\alpha)}y(x_j,t)},\label{eq07}\\
\frac{\partial^\alpha y(x_i,t)}{\partial_- x^\alpha}\cong\sum_{k=0}^M\delta_k(t)\frac{\partial^\alpha \phi_k(x_i)}{\partial_- x^\alpha}
        =\sum_{k=0}^M\delta_k(t)\sum\limits_{j=0}^M {b_{ij}^{(\alpha)}\phi_k(x_j)}\cong \sum\limits_{j=0}^M {b_{ij}^{(\alpha)}y(x_j,t)},\label{eq08}
\end{align}
\end{small}
thanks to the linearity of the fractional derivatives, namely,  Eqs. (\ref{eq07})-(\ref{eq08}) are valid 
as along as Eqs. (\ref{eq05})-(\ref{eq06}) are satisfied. The idea behind this  is referred to as DQ \cite{R30};
$a_{ij}^{(\alpha)}$, $b_{ij}^{(\alpha)}$, $i,j=0,1,\ldots,M$, are called the weighted coefficients of fractional derivatives and
will be calculated by means of typical RBFs.

\subsection{Three typical RBFs}
RBFs are the functions of the distance from their centers. They are popular as an effective tool to set up numerical algorithms
for PDEs since the superiority of potential spectral accuracy. Here, commonly used RBFs are involved, i.e.,  
\begin{itemize}
  \item Multiquadrics (MQ):  $\varphi_k(x)=\sqrt{r_k^2+\epsilon^2}$
  \item Inverse Multiquadrics (IM): $\varphi_k(x)=\frac{1}{\sqrt{r_k^2+\epsilon^2}}$
  \item Gaussians (GA):  $\varphi_k(x)=e^{-\epsilon r_k^2}$
\end{itemize}
where $r_k=|x-x_k|$, $k=0,1,\ldots,M$ and $\epsilon$ is the shape parameter.
It is worthy to note that the value $\epsilon$ should be well prescribed in computation
because it has a significant impact on the approximation power of a RBFs based method.

\subsection{Weighted coefficients for fractional derivatives}
In order to obtain the weighted coefficients of the left and right fractional derivatives,
we substitute the RBFs into Eqs. (\ref{eq05})-(\ref{eq06}) to get
\begin{align}
\frac{\partial^\alpha \varphi_k(x_i)}{\partial_+ x^\alpha}&=\sum\limits_{j=0}^M {a_{ij}^{(\alpha)}\varphi_k(x_j)},\ \ i,k=0,1,\ldots M,\label{eq11}\\
\frac{\partial^\alpha \varphi_k(x_i)}{\partial_- x^\alpha}&=\sum\limits_{j=0}^M {b_{ij}^{(\alpha)}\varphi_k(x_j)},\ \ i,k=0,1,\ldots M.\label{eq12}
\end{align}
Rewriting Eqs. (\ref{eq11})-(\ref{eq12}) in a matrix-vector form for each grid point $x_i$ yields
\begin{align}\label{eq09}
\underbrace{\left( \begin{array}{cccc}
\varphi_0(x_0)&\varphi_0(x_1) &\cdots &\varphi_0(x_M)\\
\varphi_1(x_0)&\varphi_1(x_1) &\cdots &\varphi_1(x_M)  \\
 \vdots &\vdots &\ddots&\vdots \\
\varphi_M(x_0) &\varphi_M(x_1)&\cdots &\varphi_M(x_M)
\end{array} \right)}_{\textbf{M}}
\left( \begin{array}{c}
\omega^{(\alpha)}_{i0}\\
\omega^{(\alpha)}_{i1}\\
\vdots \\
\omega^{(\alpha)}_{iM}
\end{array} \right)
=\underbrace{\left( \begin{array}{c}
\boldsymbol{D}^\alpha\varphi_0(x_i)\\
\boldsymbol{D}^\alpha\varphi_1(x_i)\\
\vdots \\
\boldsymbol{D}^\alpha\varphi_M(x_i)
\end{array} \right)}_{\boldsymbol{D}^\alpha\boldsymbol{\varphi}(x_i)},
\end{align}
where $\omega^{(\alpha)}_{ij}=a^{(\alpha)}_{ij}$ if $\boldsymbol{D}^\alpha=\frac{\partial^\alpha}{\partial_+ x^\alpha}$ whereas
$\omega_{ij}^{(\alpha)}=b^{(\alpha)}_{ij}$ if $\boldsymbol{D}^\alpha=\frac{\partial^\alpha}{\partial_- x^\alpha}$,
$i,j=0,1,\ldots,M$. $\textbf{M}$ is the interpolation matrix only related to the nodal distribution,
being nonsingular for MQ and fully positive definite for IM, GA \cite{R31}. One has
\begin{equation*}
\textbf{M}=\left( \begin{array}{cccc}
\epsilon &\sqrt{(x_1-x_0)^2+\epsilon^2} &\cdots &\sqrt{(x_M-x_0)^2+\epsilon^2}\\
\sqrt{(x_0-x_1)^2+\epsilon^2}&\epsilon &\cdots &\sqrt{(x_M-x_1)^2+\epsilon^2}\\
\vdots &\vdots &\ddots&\vdots \\
\sqrt{(x_0-x_M)^2+\epsilon^2}&\sqrt{(x_1-x_M)^2+\epsilon^2}&\cdots &\epsilon
\end{array} \right),
\end{equation*}
in particular, when MQ RBFs are adopted.

There are clearly no explicit expressions for $\frac{\partial^\alpha\varphi_k(x)}{\partial_+ x^\alpha}$,
$\frac{\partial^\alpha\varphi_k(x)}{\partial_- x^\alpha}$; fortunately, they can be approximated by numerical quadrature rules.
Taking the transforms $\xi=x-\frac{(x-a)(1+\zeta)}{2}$, $\xi=x+\frac{(b-x)(1+\zeta)}{2}$ of variables, respectively, reaches to
\begin{align*}
 \frac{\partial^\alpha \varphi_k(x_i)}{\partial_+ x^\alpha}&=\frac{1}{\Gamma(2-\alpha)}\bigg(\frac{x_i-a}{2}\bigg)^{2-\alpha}
    \int^1_{-1}(1+\zeta)^{1-\alpha}\varphi^{''}_k\bigg(x_i-\frac{(x_i-a)(1+\zeta)}{2}\bigg)d\zeta,\\
 \frac{\partial^\alpha \varphi_k(x_i)}{\partial_- x^\alpha}&=\frac{1}{\Gamma(2-\alpha)}\bigg(\frac{b-x_i}{2}\bigg)^{2-\alpha}
    \int^1_{-1}(1+\zeta)^{1-\alpha}\varphi^{''}_k\bigg(x_i+\frac{(b-x_i)(1+\zeta)}{2}\bigg)d\zeta.
\end{align*}
A close examination reveals that both of the two formulas are the special cases of the following weakly singular integral, i.e., 
\begin{equation*}
\int^1_{-1}(1-\zeta)^\lambda(1+\zeta)^{\mu}f(\zeta)d\zeta,\quad \lambda,\ \mu>-1,
\end{equation*}
with $\lambda=0$, $\mu=1-\alpha$ that can be handled by Gauss-Jacobi quadrature rules. $a_{ij}^{(\alpha)}$, $b_{ij}^{(\alpha)}$
are then determined by solving Eqs. (\ref{eq09}) for each  $x_i$ and the fractional derivatives are
removed from a fractional PDE by using Eqs. (\ref{eq07})-(\ref{eq08}) as replacements, thus we obtain the required solution by solving a ODS instead.

\section{A time-stepping RBF-DQM for fractional PDEs}\label{s4}
In this section, a RBF-DQM of fully discretization is derived for the space-fractional PDEs via the above direct approximations for fractional derivatives.
Define a lattice on $[0,T]$ with equally spaced points $t_n=n\tau$, $\tau=T/N$, $N\in\mathbb{Z}^+$.
On substituting the weighted sums (\ref{eq07})-(\ref{eq08})
into Eq. (\ref{eq01}), we have
\begin{align*}
\frac{\partial y(x_i,t)}{\partial t}-\kappa(x_i)\sum\limits_{j=0}^M {a_{ij}^{(\alpha)}y(x_j,t)}
 -\upsilon(x_i)\sum\limits_{j=0}^M {b_{ij}^{(\alpha)}y(x_j,t)}=f(x_i,t), \ \ i=0,1,\cdots,M,
\end{align*}
actually being a first-order ODS. 
Also, denote $t_{n-1/2}=t_n-\frac{\tau}{2}$, $y^{n}_i=y(x_i,t_{n})$, $f^{n-1/2}_i=f(x_i,t_{n-1/2})$ for brevity.
Imposing the associated constraints (\ref{eq02})-(\ref{eq03}) and rewriting the ODS in matrix-vector form,
a time-stepping RBF-DQM is then derived by introducing a Crank-Nicolson scheme in time, given as
\begin{equation}\label{eq10}
    \bigg(\textbf{I}-\tau\frac{\boldsymbol{\kappa}\textbf{A}+\boldsymbol{\upsilon}\textbf{B}}{2}\bigg)\textbf{Y}^{n}=
    \bigg(\textbf{I}+\tau\frac{\boldsymbol{\kappa}\textbf{A}+\boldsymbol{\upsilon}\textbf{B}}{2}\bigg)\textbf{Y}^{n-1}
        +\tau\textbf{H}^{n-1/2}, 
\end{equation}
where \textbf{I} is an identity matrix,
$\textbf{Y}^{n}=[y^n_1,y^n_2,\cdots,y^n_{M-1}]^T$,
$\boldsymbol{\kappa}=\textrm{diag}(\kappa_1,\kappa_2,\cdots,\kappa_{M-1})$,
$\boldsymbol{\upsilon}=\textrm{diag}(\upsilon_1,\upsilon_2,\cdots,\upsilon_{M-1})$,
and $\textbf{A}$, $\textbf{B}$, $\textbf{H}^{n-1/2}$ are as follows
\begin{align*}
&\textbf{A}=\left( \begin{array}{cccc}
a^{(\alpha)}_{11}&a^{(\alpha)}_{12} &\cdots &a^{(\alpha)}_{1,M-1}\\
a^{(\alpha)}_{21}&a^{(\alpha)}_{22} &\cdots &a^{(\alpha)}_{2,M-1}  \\
 \vdots &\vdots &\ddots&\vdots \\
a^{(\alpha)}_{M-1,1}&a^{(\alpha)}_{M-1,2}&\cdots &a^{(\alpha)}_{M-1,M-1}
\end{array} \right),
\textbf{B}=\left( \begin{array}{cccc}
b^{(\alpha)}_{11}&b^{(\alpha)}_{12} &\cdots &b^{(\alpha)}_{1,M-1}\\
b^{(\alpha)}_{21}&b^{(\alpha)}_{22} &\cdots &b^{(\alpha)}_{2,M-1}  \\
 \vdots &\vdots &\ddots&\vdots \\
b^{(\alpha)}_{M-1,1}&b^{(\alpha)}_{M-1,2}&\cdots &b^{(\alpha)}_{M-1,M-1}
\end{array} \right),\\
&\textbf{H}^{n-1/2}=\left( \begin{array}{c}
f^{n-1/2}_1\\
f^{n-1/2}_2\\
 \vdots \\
f^{n-1/2}_{M-1}
\end{array} \right)
+\frac{g_0^n+g_0^{n-1}}{2}\left( \begin{array}{c}
\omega^{(\alpha)}_{1}\\
\omega^{(\alpha)}_{2}\\
 \vdots \\
\omega^{(\alpha)}_{M-1}
\end{array} \right)
+\frac{g_M^n+g_M^{n-1}}{2}\left( \begin{array}{c}
\tilde{\omega}^{(\alpha)}_{1}\\
\tilde{\omega}^{(\alpha)}_{2}\\
 \vdots \\
\tilde{\omega}^{(\alpha)}_{M-1}
\end{array} \right),
\end{align*}
with $\omega^{(\alpha)}_{i}=\kappa_ia^{(\alpha)}_{i0}+\upsilon_ib^{(\alpha)}_{i0}$,
$\tilde{\omega}^{(\alpha)}_{i}=\kappa_ia^{(\alpha)}_{iM}+\upsilon_ib^{(\alpha)}_{iM}$.
We perform the procedures on the nodal distribution
$x_i=0.5(1-\cos\frac{i\pi}{M})\ell+a$, $\ell=b-a$, $i=0,1,\cdots,M$.
An detailed implementation of RBF-DQM is summarized in the following flowchart
\begin{itemize}
  \item Input $\alpha$, $\epsilon$, $M$, $N$, and allocate $\{t_n\}_{n=0}^{N}$, $\{x_i\}_{i=0}^{M}$.
  \item Form \textbf{M}, compute $\boldsymbol{D}^\alpha\boldsymbol{\varphi}(x_i)$ by Gauss-Jacobi quadrature rules,
         and solve Eqs. (\ref{eq09}) for each $x_i$ so that the weighted coefficients are found.
  \item Do a loop from $n=1$ to $N$ to solve  Eqs. (\ref{eq10}) for each $t_n$ by forming \textbf{A}, \textbf{B} first and output the desirable approximation $\textbf{Y}^{n}$ at each time step.
\end{itemize}

\section{Illustrative examples}\label{s5}

In this part, the proposed methods, termed as MQ-DQM, IM-DQM, and GA-DQM hereinafter, are studied
on three numerical examples. The shape parameter should be adjusted with the grid number $M$,
so we select $\epsilon=1.25\ell/(M+1)^{0.5}$ for MQ, $\epsilon=2/(M+1)^{0.5}$
for IM, and $\epsilon=1.05(M+1)$ for GA, tentatively, by references to \cite{R62,R32,R61}. The numerical errors are
all defined in $l_\infty$-norm and the fractional derivatives are computed with 16 quadrature points and weights,
whose values corresponding to $\alpha=1.5$ are given as reference in \ref{appd}. \\ 


\noindent
\textbf{Example 5.1.}
Approximate ${_0}D_x^{1.1}\sin(x)$ by Eq. (\ref{eq07}) on $[0,1]$ with the above $\epsilon$, whose explicit
expression is $-x^{1.9}{_1}F_2(1;1.45,1.95;-0.25x^2)/\Gamma(2.9)$, where ${_1}F_2$ is
the \emph{hypergeometric function}. The numerical results are tabulated in Table \ref{tab1}.
As observed, the approximation improves as $M$ increases, which implies that the DQ approximations
for the fractional derivatives are valid. Moreover, under the given $\epsilon$,
GA-DQM seems to be more efficient than MQ-DQM and IM-DQM.\\


\begin{table*}[!htb]
\centering
\caption{The numerical results when $\alpha=1.1$ for Example 5.1} \label{tab1}
\begin{tabular}{cccc}
\hline $M$ &MQ-DQM &IM-DQM &GA-DQM \\
\hline  5  & 6.4167e-02  &5.7172e-02  &1.9069e-01  \\
        10 & 6.3543e-03  &8.5795e-03  &2.2488e-02  \\
        15 & 1.3676e-03  &7.2670e-04  &8.5354e-04  \\
        20 & 4.2484e-04  &4.0834e-04  &1.1827e-04  \\
        25 & 2.2093e-04  &8.2798e-05  &4.6350e-06  \\
\hline
\end{tabular}
\end{table*}

\noindent
\textbf{Example 5.2.} Let $\kappa(x)=\upsilon(x)=1$, and $y(x,t)=t^3x^2(1-x)^2$;
we solve Eqs. (\ref{eq01})-(\ref{eq03}) on $[0,1]$ with zero initial-boundary conditions.
Table \ref{tab2} displays the numerical results at $t=0.5$ when $\tau=2.0\times10^{-4}$ and $\alpha=1.8$.
From the table, we find that sufficiently small errors can be achieved by MQ-DQM, IM-DQM, and GA-DQM even if a
few spatial grid points are utilized and all the methods are obviously
convergent by taking their own $\epsilon$, respectively. \\

\noindent
\textbf{Example 5.3.} Let $\kappa(x)=\frac{x^\alpha\Gamma(5-\alpha)}{24}$, $\upsilon(x)=0$, and $y(x,t)=e^{-t}x^4$;
we solve Eqs. (\ref{eq01})-(\ref{eq03}) on $[0,1]$ with $\psi(x)=x^4$, $g_1(t)=0$, and $g_2(t)=e^{-t}$.
The numerical errors of SAM \cite{R17} and our methods at $t=1$ are reported in Table \ref{tab3},
when $\tau=1/M$ and $\alpha=1.5$. It is observed from the table that MQ-DQM, IM-DQM, and GA-DQM
outperform SAM in term of computational accuracy. \\

\begin{table*}[!htb]
\centering
\caption{The numerical results at $t=0.5$ when $\tau=2.0\times10^{-4}$, $\alpha=1.8$ for Example 5.2} \label{tab2}
\begin{tabular}{ccccc}
\hline  $M$&  MQ-DQM &IM-DQM & GA-DQM \\
\hline  5  &  2.2567e-04  &  2.0463e-04  &2.2607e-04  \\
        10 &  1.5291e-05  &  1.2391e-05  &7.3949e-06  \\
        20 &  4.1822e-07  &  1.9394e-07  &1.4895e-08  \\
        25 &  7.6704e-08  &  2.9039e-08  &2.0757e-09  \\
\hline
\end{tabular}
\end{table*}

\begin{table*}[!htb]
\centering
\caption{A comparison of SAM and RBF-DQM  at $t=1$ when $\tau=1/M$ and $\alpha=1.5$.} \label{tab3}
\begin{tabular}{ccccc}
\hline  $M$ & SAM \cite{R17} & MQ-DQM &IM-DQM & GA-DQM \\
\hline  15 & 7.660e-04 & 1.5903e-04  &  1.4862e-04  &1.4127e-04  \\
        20 & 4.493e-04 & 7.9355e-05  &  7.5052e-05  &7.2001e-05  \\
        25 & 2.929e-04 & 4.6347e-05  &  4.5794e-05  &4.5247e-05  \\
        30 & 2.067e-04 & 2.9290e-05  &  3.0308e-05  &3.1650e-05  \\
\hline
\end{tabular}
\end{table*}

\begin{remark}
When $M$ is fixed, the value $\epsilon$ is crucial to the accuracy of a RBFs based method, so is RBF-DQM.
 A general trade-off principle demanding attention is that one can adjust $\epsilon$ to decrease the approximate errors, 
but need to pay for this by increasing the condition number of the interpolation matrix which may cause an algorithm to be instable \cite{R33},
so a good $\epsilon$ that balances both the accuracy and stability is anticipated in practice.
Nevertheless, how to select an optimal value is technical and is being an issue deserving to investigate.
\end{remark}

\section{Conclusion}
In this research, an efficient DQ method is proposed for the space-fractional PDEs of Caputo type based on
commonly used RBFs as test functions, which enjoys some properties such as high accuracy, flexibility,
truly meshless, and the simplicity in implementation. Its codes are tested on three benchmark examples and the
outcomes manifest that it is capable of dealing with these problems if the free parameters $\epsilon$ are well prepared.
Due to its insensitivity to dimensional change, our method has potential advantages over traditional methods
in finding the approximate solutions to the high-dimensional fractional equations. \\


\noindent
\textbf{Acknowledgement}:
This research was supported by National Natural Science Foundations of China (No.11471262 and 11501450).

\begin{appendix}
\section{Quadrature points and weights when $\alpha=1.5$ }\label{appd}
\begin{table*}[!htb]
\small
\centering
\begin{tabular}{rrrr}
\hline   \multicolumn{2}{l}{quadrature points}&\multicolumn{2}{l}{weights} \\
\hline   -0.995332738871603 & 0.072181310006040  & 0.273056604980456  &0.186173360311749 \\
         -0.958255704632838 & 0.262409004111262  & 0.270507150048226  &0.165983681785832 \\
         -0.885483116532829 & 0.442850686515520  & 0.265432043865656  &0.144244294127118 \\
         -0.779726461638614 & 0.606783228547583  & 0.257878671571417  &0.121158215878927 \\
         -0.644926204206171 & 0.748098741650606  & 0.247917557394430  &0.096941123533022 \\
         -0.486104964648995 & 0.861532419612587  & 0.235641706298296  &0.071819606274043 \\
         -0.309180377381835 & 0.942860103381098  & 0.221165735833726  &0.046030939495720 \\
         -0.120744600181947 & 0.989070420301884  & 0.204624806418771  &0.019851626928800 \\
\hline
\end{tabular}
\end{table*}
\end{appendix}

\bibliographystyle{model1b-num-names}
\bibliography{mybib}

\end{document}